\documentclass[11pt]{amsart}
\usepackage[T1]{fontenc}
\usepackage[ansinew]{inputenc}
\usepackage{amssymb}
\newtheorem{thm}{Theorem}[section]

\newtheorem{cor}[thm]{Corollary}
\newtheorem{lem}[thm]{Lemma}
\newtheorem{prop}[thm]{Proposition}
\newtheorem{defn}[thm]{Definition}
\newtheorem{rem}{Remark}[section]
\newcommand{\eps}{\varepsilon}

\newcommand{\by}{\bar{y}}

\newcommand{\ut}{\underline{t}}
\newcommand{\bt}{\bar{t}}
\begin{document}
\title{Controllability of 2D Euler and Navier-Stokes Equations
by Forcing $4$ Modes}%
\author{Andrey A. Agrachev$^1$ \and Andrey
V. Sarychev$^2$}%
\address{$^1$International School for Advanced Studies (SISSA),
Trieste, Italy \& V.A.Steklov Mathematical Institute, Moscow, Russia
\newline  $^2$ DiMaD,
University of Florence, Italy}%
\email{agrachev@sissa.it,asarychev@unifi.it}%

\date{}%
%\dedicatory{}%
%\commby{}%
\begin{abstract}
We  study  controllability issues for the 2D Euler and
Navier-Stokes (NS) systems under periodic boundary conditions.
These  systems describe motion  of homogeneous ideal or viscous
incompressible fluid on a two-dimensional torus $\mathbb{T}^2$. We
assume the system to be controlled by a degenerate forcing applied
to fixed number of modes.

In our previous work \cite{ASpb,AS43,ASDAN} we studied global
controllability by means of degenerate forcing for Navier-Stokes
(NS) systems with nonvanishing viscosity ($\nu
>0$). Methods of differential geometric/Lie algebraic control
theory have been used for that study. In \cite{ASpb} criteria for
global controllability of finite-dimensional Galerkin
approximations of 2D and 3D NS systems have been established. It
is almost immediate to see that these criteria are also valid for
the Galerkin approximations of the Euler systems; in fact the Lie
brackets involved into the corresponding Lie rank controllability
condition do not depend on the viscous term. In \cite{AS43,ASDAN}
we established a much more intricate sufficient criteria for
global controllability in finite-dimensional observed component
and for $L_2$-approximate controllability for 2D NS system . The
justification of these criteria was based on a Lyapunov-Schmidt
reduction to a finite-dimensional system. Possibility of such a
reduction rested upon the dissipativity of NS system, and hence
the previous approach can not be adapted for Euler system.

 In the present contribution we improve and extend the controllability
results in several aspects: 1) we obtain a stronger sufficient
condition for controllability of 2D NS system in an observed
component  and for $L_2$-approximate controllability; 2) we prove
that these criteria are valid for the case of ideal incompressible
fluid ($\nu=0)$; 3) we study solid controllability in projection
on any finite-dimensional subspace and establish a sufficient
criterion for such controllability.
\end{abstract}
% ----------------------------------------------------------------
\maketitle
% ----------------------------------------------------------------
%%%%%%%%%%%%%%%%%%%%%%%%%%%%%%%%%%%%%%%%%%%%%%%%%%%%%%%%%%%%%%%%%%%%

\vspace*{2mm}


\section{Introduction}
\label{intro} \markboth{A.A.Agrachev,
A.V.Sarychev}{Controllability of 2D Euler and NS Equations by $4$
Modes Forcing}

The present paper extends our work  started in
 \cite{ASpb,AS43,ASDAN}
 on studying controllability of  $2$- and $3$-
dimensional Navier-Stokes equations (2D and 3D NS systems) under
periodic boundary conditions. The characteristic feature of our
problem setting is  a choice of control functions; we are going to
control the 2D NS/Euler system  by means of {\it degenerate}
forcing. The corresponding equations are
\begin{eqnarray}\label{eul1}
\partial u/\partial t + (u \cdot \nabla)u + \nabla p = \nu \Delta u +
F(t,x),  \\
\nabla \cdot u=0.\label{ns2}
\end{eqnarray}
The words "degenerate forcing" mean that $F(t,x)$ is a "low-order"
trigonometric polynomial with respect to $x$, i.e. a sum of a
"small number" of harmonics:
$$F(t,x)=\sum_{k \in {\mathcal
K}^1}v_k(t)e^{ik \cdot x}, \ {\mathcal K}^1 \ \mbox{is finite}.$$
The word "control" means that the components $v_k(t), \ t \in
[0,T]$ of the forcing can be chosen freely among measurable
essentially bounded functions. In fact to achieve controllability
piecewise-constant controls suffice.

In \cite{ASpb,AS43,ASDAN} we derived sufficient controllability
criteria  for Galerkin approximations of 2D and 3D NS systems. For
the 2D NS system we established sufficient criteria for so called,
controllability in finite-dimensional observed component and for
$L_2$-approximate controllability. The corresponding definitions
can be found in the Section~\ref{dset}.

Now we consider both  cases of viscous ($\nu
>0$) and  {\em ideal} ($\nu=0$) incompressible fluid
simultaneously. To establish a possibility to propagate the action
of small dimensional control to (a finite number of) higher modes
we use the technique of Lie extensions developed in the scope of
geometric control theory (see \cite{ASkv,Ju}). For
finite-dimensional Galerkin approximations of 2D and 3D Euler
systems ($\nu =0$) the controllability criteria turn out to be
the same as for 2D and 3D NS systems ($\nu
>0$)  (see \cite{ASpb,AS43,ASDAN}).
This is due to the fact that these controllability criteria are
of "purely nonlinear" nature; they are completely determined by
the nonlinear term of the Euler system.

Tools of Geometric Control Theory are not  yet adapted too much to
infinite-dimensional case. For dealing with inifinite-dimensional
dynamics we used in \cite{ASpb,AS43,ASDAN} a Lya\-pu\-nov-Schmidt
reduction to a finite-dimensional system. The possibility of such
a reduction rested upon dissipativity of the NS system, which is
not anymore present when one deals with Euler system.

In the present paper we abandon the Lya\-pu\-nov-Schmidt reduction
and instead refine the tools of geometric control in order to deal
with viscous and nonviscous case at the same time. This refinement
also allow us to improve the sufficient criterion of
controllability in observed component for 2D NS/Euler system. The
criterion, formulated in terms of so-called 'saturating property'
of the set of controlled forcing modes, is stronger then the one
established in \cite{ASpb,AS43,ASDAN}. Analysis of the saturation
property in \cite{HM} showed that a generic symmetric set of $4$
controlled modes suffices for achieving controllability.

For a saturating set we manage to prove $L_2$-approximate
controllability for 2D NS/Euler system. We also study
controllability in finite-dimensional projections. The latter
property means that the attainable set of 2D NS/Euler system is
projected surjectively onto any finite-dimensional subspace of
$H_2$.

There has been an extensive study of controllability of the
Navier-Stokes and Euler equations in particular by means of
boundary control. There are various results on exact local
controllability of 2D and 3D Navier-Stokes equations obtained by
A.Fursikov, O.Imanuilov, global exact controllability for 2D Euler
equation obtained by J.-M. Coron, global exact controllability for
2D Navier-Stokes equation by A.Fursikov and J.-M. Coron. We refer
the readers to the book \cite{Fu} and to the surveys \cite{FI} and
\cite{Co} for further references.

Our problem setting differs from the above results by the class of
{\em degenerate distributed controls}  which is involved. In
closer relation to our work is a publication of M.Romito
(\cite{Ro}) who provided a criterion for controllability of
Galerkin approximations of 3D NS systems.  J.C.Mattingly and
E.Pardoux adapted  (\cite{MP})) the controllability result from
\cite{AS43} for  studying properties of the solutions of
stochastically forced 2D NS systems. M.Hairer and J.C.Mattingly
have applied (\cite{HM}) the controllability results to studying
ergodicity of 2D Navier-Stokes equation under degenerate
stochastic forcing.

The structure of our paper is as follows. Section~\ref{prel}
contains a necessary minimum of standard preliminary material on
2D Euler and NS systems. The problem setting in the
Section~\ref{dset} is succeeded by the formulation of the main
results in the Section~\ref{mres}. These results include
sufficient criteria for controllability in observed component, for
solid controllability in finite-dimensional projection, and for
$L_2$-approximate controllability for both 2D NS and 2D Euler
systems.

The rest of the paper is devoted to the proofs of these results.
Among the tools involved are some results on equiboundedness and
continuous dependence of solutions of 2D NS/Euler  systems on
relaxed forcings. Being interesting for their own sake these
results are formulated in Section~\ref{cont} and are proved in
Appendix. In  Section~\ref{cr} we accomplish the proof of (solid)
controllability in observed component. The construction,
introduced in this proof, is crucial for the proof of
controllability in a finite-dimensional projection accomplished in
Section~\ref{proc}. Proof of $L_2$-approximate controllability is
similar; the readers can either complete it by themselves or
consult \cite{AS43}.

%%%%%%%%%%%%%%%%%%%%%%%%%%%%%%%%%%%%%%%%%%%%%%%%%%%%%%%%%%%%%%%%%%%%%

\section{Preliminaries on 2D NS/Euler System: vorticity, spectral method, Galerkin app\-ro\-xi\-mations}
\label{prel}

We consider the 2D NS/Euler system (\ref{eul1})-(\ref{ns2}). The
boundary conditions are assumed to be periodic, i.e. one may
assume the velocity field $u$ to be defined on the
$2$-dimensional torus $\mathbb{T}^2$. Besides we assume
\begin{equation}\label{zerav}
\int_{\mathbb{T}^2}udx=0.
\end{equation}

Let us introduce the vorticity $w=\nabla^\perp \cdot u=\partial
u_2/\partial x_1 - \partial u_1/\partial x_2$ of $u$. Applying the
operator $\nabla^\perp$ to the equation (\ref{eul1}) we arrive to
the  equation:
\begin{equation}\label{nsw}
\partial w/\partial t + (u \cdot \nabla)w -\nu \Delta w =
v(t,x),
\end{equation}
where $v(t,x)=\nabla^\perp \cdot F(t,x)$.

Notice that: i) $\nabla^\perp \cdot \nabla p=0$, ii)
$\nabla^\perp$ and $\Delta$ commute as linear differential
operators in $x$ with constant coefficients;
 $$\mbox{iii)} \  \nabla^\perp
\cdot (u \cdot \nabla)u= (u \cdot \nabla)(\nabla^\perp \cdot
u)+(\nabla^\perp \cdot u)(\nabla \cdot u)=(u \cdot \nabla)w,$$ for
all $u$ satisfying (\ref{ns2}).

It is known that $u$, which satisfies the relations (\ref{ns2})
and (\ref{zerav}), can be recovered in a unique way from $w$. From
now on we will deal with the equation (\ref{nsw}).

A natural and standard (see \cite{BaV,CF}) way to view the NS
systems is to represent them  as evolution equations in Hilbert
spaces.

Consider Sobolev spaces $H^\ell(\mathbb{T}^s)$ with the scalar
product defined as
$$\langle u,u'\rangle_\ell  =\sum_{\alpha \leq
\ell}\int_{\mathbb{T}^s}(\partial^\alpha u/\partial x^\alpha)
(\partial^\alpha u'/\partial x^\alpha) dx;$$ the norm $\| \cdot
\|_\ell$ is defined by virtue of this scalar product. Denote by
$H_\ell$ the closures of $\{u \in C^\infty(\mathbb{T}^s), \nabla
\cdot u=0\}$ in the norms $\|\cdot\|_\ell$ in the respective
spaces $H^\ell(\mathbb{T}^s), \ \ell \geq 0$.    The norms in
$H_\ell$ will be denoted again by $\|\cdot\|_\ell$.   It will be
convenient for us to redefine the norm of $H_1$ by putting
$\|u\|^2_1=\langle -\Delta u, u\rangle$, and the norm of $H_2$ by
putting $\|u\|^2_2=\langle -\Delta u, - \Delta u\rangle$.

Results on global existence and uniqueness of weak and classical
solutions of NS systems in bounded domains can be found in
\cite{CF,BaV,Lad}. The proofs for the the non-viscous (Euler) case
are more delicate. W.Wolibner's existence and uniqueness theorem
is presented in \cite{Ka}. Formulation in \cite{EMar} allows for
asserting global existence and uniqueness of trajectories $t
\mapsto u_t$ (respectively $t \mapsto w_t$ for the vorticity) of
the 2D Euler systems in any Sobolev space $H_s$ with $s>2$
(respectively with $s>1$ for the vorticity), provided that the
initial data belongs to these spaces.

Let us consider now the basis of eigenfunctions $\{e^{ik \cdot
x}\}$ of the Laplacian on $\mathbb{T}^2$ and take the Fourier
expansion of the vorticity  $w(t,x)=\sum_k q_k(t)e^{ik \cdot x}$
and control $v(t,x)=\sum_k v_k(t)e^{ik \cdot x}$. As far as $w$
and $f$ are real-valued, we have $\bar{w}_n=w_{-n}, \
\bar{v}_n=v_{-n}$. We assume $v_0=0$; by (\ref{zerav}) $q_0=0$.

 Evidently $\partial w/\partial t = \sum_k \dot{q}_k(t)e^{ik\cdot
 x}$.
To compute $(u \cdot \nabla)w$ we write the equalities
$$\nabla^\perp \cdot u=w, \ \nabla \cdot u=0 \Leftrightarrow  -\partial_2 u_1+\partial_1 u_2=w, \
\partial_1 u_1+\partial_2 u_2=0.$$
From these latter we conclude by a standard reasoning that
$$u_1=\sum_k q_k(t)(ik_2/|k|^2)e^{ik \cdot x}, \ u_2=-\sum_k
q_k(t)(ik_1/|k|^2)e^{ik \cdot x},$$ and
$$(u \cdot \nabla)w=\sum_{m+n=k}(m \wedge n)|m|^{-2}q_m q_n,$$
where $m \wedge n=m_1n_2-m_2n_1$ is the external product of
 $m=(m_1,m_2), n=(n_1,n_2)$.

Now the 2D NS/Euler system can be written as  an
(infinite-dimensional) system of ODE for $q_k$:
\begin{equation}\label{1mais}
\dot{q}_k=\sum_{m+n=k}(m \wedge n)|m|^{-2}q_m q_n-\nu
|k|^2q_k+v_k, \ k,m,n \in \mathbb{Z}^2.
\end{equation}

Observe that the product $q_mq_n$ enters the sum $\sum_{m+n=k}(m
\wedge n)|m|^{-2}q_m q_n$ twice with (a priori) different
coefficients. Therefore this sum can be rearranged
\begin{equation}\label{mleqn}
\sum_{m+n=k}(m \wedge n)|m|^{-2}q_m q_n=\sum_{m+n=k, |m|<|n|}(m
\wedge n)(|m|^{-2}-|n|^{-2})q_m q_n.
\end{equation}
From the last representation we conclude that $q_mq_n$ does not
appear in the equations whenever $|m|=|n|$.

Consider any finite subset $ \mathcal{G} \subset \mathbb{Z}^2$ and
introduce the Galerkin $\mathcal{G}$-ap\-pro\-xi\-ma\-tion of the
system (\ref{nsw}) or of the system (\ref{1mais}) by projecting
this system  onto the linear space spanned by the harmonics
$e^{ik\cdot x}$ with $k \in \mathcal{G}$. The result is a
finite-dimensional system of ODE or a control system
\begin{equation}\label{1gal}
\dot{q}_k=\sum_{m+n=k}(m \wedge n)|m|^{-2}q_m q_n-\nu
|k|^2q_k+v_k, \ k,m,n \in \mathcal{G}.
\end{equation}

\section{2D NS/Euler system controlled by degenerate forcing.
Problem setting} \label{dset}

We study the case where the 2D NS/Euler  system is forced by a
trigonometric polynomial: $v(t,x)=\sum_{k \in {\mathcal
K}^1}v_ke^{i k \cdot x}$, where ${\mathcal K}^1$ is a finite set.
Such forcing is called {\em degenerate}. As we said $v_k(\cdot)$
with $k \in {\mathcal K}^1$ are controls at our disposal; they are
arbitrary measurable essentially bounded functions.

Let us introduce a finite set of {\em observed} modes indexed by
$k \in \mathcal{K}^{obs} \subset \mathbb{Z}^2$. The observed modes
are reunited in so-called observed component.

We assume $\mathcal{K}^{obs} \supseteq \mathcal{K}^1$. As we will
see, nontrivial controllability issues arise only if
$\mathcal{K}^1$ is  a proper subset of $\mathcal{K}^{obs}$.
  We identify the space of observed modes with
$\mathbb{R}^N$ and denote by $\Pi^{obs}$ the operator of
projection of solutions onto the space of observed modes.

We will represent the controlled 2D NS/Euler equation in the
following splitted (controlled -observed -unobserved components)
form:
\begin{eqnarray}\label{2con1}
\dot{q}_k=\sum_{m+n=k}(m \wedge n)|m|^{-2}q_m q_n-\nu |k|^2q_k+v_k, \ k \in \mathcal{K}^1,  \\
\label{2cobs} \dot{q}_k=\sum_{m+n=k}(m \wedge n)|m|^{-2}q_m q_n-\nu |k|^2q_k, \ k \in \mathcal{K}^{obs} \setminus \mathcal{K}^1, \\
\label{2unobs} \dot{Q}=\nu \Delta Q+B(q,Q).
\end{eqnarray}

In the latter equation $B(q,Q)$ stays for the projection of the
 nonlinear term of the Euler system onto the space of unobserved
modes.

Galerkin $\mathcal{K}^{obs}$-approximation of the 2D NS/Euler
system consists of the equations (\ref{2con1})-(\ref{2cobs}) under
an additional condition $m,n \in \mathcal{K}^{obs}$ for the
summation indices.

\begin{defn}
Galerkin $\mathcal{K}^{obs}$-approximation of 2D  NS/Euler systems
is time-$T$ globally controllable if for any two points
$\tilde{q}, \hat{q}$  in $\mathbb{R}^N$,  there exists a control
which steers in time $T$ this Galerkin approximation from
$\tilde{q}$ to $\hat{q}. \ \Box$
\end{defn}

\begin{defn}{\it (controllability in observed component)}
\label{gcop}
 2D NS/Euler  system is  time-$T$ globally
controllable in observed $N$-dimensional component  if for any
$\tilde{\varphi} \in H_2$ and any $\hat{q} \in \mathbb{R}^N$ there
exists a control which steers the system in time $T$ from
$\tilde{\varphi}$ to some $\hat{\varphi} \in
(\Pi^{obs})^{-1}(\hat{q}). \ \Box$
\end{defn}

In other words the 2D NS/Euler system is  globally controllable in
observed component if its time-$T$ attainable set (from each
point) is projected by $\Pi^{obs}$ onto the whole coordinate
subspace spanned by the observed modes. Notice that the observed
dynamics (\ref{2con1})-(\ref{2cobs}) is affected by
in\-finite-di\-men\-sional dynamics (\ref{2unobs}).

We generalize the previous definition.

\begin{defn}{\it (controllability in finite-dimensional
projection)} \label{coproj}
 Let $\mathcal L$ be a finite-dimensional subspace of $H_2(\mathbb{T}^2)$ and $\Pi^{\mathcal L}$ be
 $L_2$-orthogonal projection of $H_2(\mathbb{T}^2)$ onto $\mathcal L$.
  The 2D NS/Euler  system is  time-$T$ globally
controllable in projection on $\mathcal L$ if for any
$\tilde{\varphi} \in H_2(\mathbb{T}^2)$ and for any
 $\hat{q} \in \mathbb{R}^N$
there exists a control which steers the  system in time $T$ from
$\tilde{\varphi}$ to some $\hat{\varphi} \in (\Pi^{\mathcal
L})^{-1}(\hat{q}). \ \Box$
\end{defn}

\begin{rem}
Controllability in observed component amounts to  controllability
in finite-dimensional projection on a {\it coordinate} subspace
${\mathcal L}. \ \Box$.
\end{rem}

\begin{defn}{\it ($L_2$-approximate controllability)}
\label{dapc}
 The 2D NS/Euler system is  time-$T$
$L_2$-appro\-xi\-mate\-ly controllable, if for any two points
$\tilde{\varphi}, \hat{\varphi} \in H_2$ and for any $\eps >0$
there exists a control which steers the  system in time $T$ from
$\tilde{\varphi}$ to the $\eps$-neighborhood of $\hat{\varphi}$ in
$L_2$-norm. $\Box$
\end{defn}

Let us introduce some useful terminology.

\begin{defn}\label{itomap}
 Fix initial condition $\tilde{\varphi} \in H_2(\mathbb{T}^2)$
 for trajectories of the controlled 2D NS/Euler system.

The correspondence between the controlled forcing $v(\cdot) \in
L_\infty\left([0,T];\mathbb{R}^d\right)$ and the corresponding
trajectory (solution) $w_t$ of the  system is established by
forcing/trajectory map ($\mathcal{F/T}$-map).

The correspondence between the controlled forcing $v(\cdot)$ and
the observed component $q(t)=\Pi^{obs}w_t$ (an
$\mathbb{R}^N$-valued function) of the corresponding trajectory is
established by forcing/observation map ($\mathcal{F/O}$-map).

If NS/Euler system is considered on an interval $[0,T]\ (T <
+\infty)$, then the map $\mathcal{F/T}_T : v(\cdot) \mapsto w_T$
is called end-point map; the map $\Pi^{obs} \circ \mathcal{F/T}_T$
is  called end-point component map, the composition $\Pi^{\mathcal
L} \circ \mathcal{F/T}_T$ is called $\mathcal L$-projected
end-point map. $\Box$
\end{defn}

\begin{rem}
In the terminology of control theory the first  two maps would be
called input/trajectory and input/output maps  correspondingly.
$\Box$
\end{rem}

\begin{rem}
  Evidently time-$T$ controllability of the NS/Euler system in
observed component or in finite-dimensional projection  is the
same as surjectiveness of the corresponding end-point maps. $\Box$
\end{rem}

Invoking these maps we will introduce a stronger notion of {\it
solid controllability}.

\begin{defn}
Let $\Phi: M^1 \mapsto M^2$ be a continuous map between two metric
spaces, and $S \subseteq M^2$ be any subset. We say that $\Phi$
covers $S$ solidly, if $S \subseteq \Phi(M^1)$ and this inclusion
is stable with respect to $C^0$-small perturbations of $\Phi$,
i.e. for some $C^0$-neighborhood $\Omega$ of $\Phi$ and for each
map $\Psi \in \Omega$, there holds: $S \subseteq \Psi(M^1). \
\Box$
\end{defn}

In what follows $M^2$ will be finite-dimensional vector space.

\begin{defn}{\it (solid controllability in finite-dimensional
projection)} \label{socoproj}
   The 2D NS/Euler  system is  time-$T$ solidly globally
controllable in projection on finite-dimensional subspace
$\mathcal L \subset H_2$, if for any bounded set $S$  in $\mathcal
L$ there exists a set of controls $B_S$ such that
$\left(\Pi^{\mathcal L} \circ \mathcal{F/T}_T\right)(B_S)$ covers
$S$ solidly. $\Box$
\end{defn}

\subsection{Problem setting}

We address the following questions.

{\bf Question 1.} Under what conditions the 2D NS/Euler system
(\ref{2con1})-(\ref{2cobs})-(\ref{2unobs}) is globally
controllable in  observed component? $\Box$

{\bf Question 2.} Under what conditions the 2D NS/Euler system
(\ref{2con1})-(\ref{2cobs})-(\ref{2unobs}) is solidly
controllable in any finite-dimensional projection? $\Box$

{\bf Question 3.} Under what conditions the 2D NS/Euler system is
$L_2$-approximately controllable? $\Box$

 In  \cite{AS43,ASDAN} we have
answered the Questions 1,3 for the 2D NS system.  In the present
contribution we improve the previous results (provide sufficient
controllability conditions under weaker hypothesi), extend them
onto the case of ideal fluid (2D Euler equation) and answer the
Question 2 for 2D NS and Euler systems.

\section{Main results for $2D$ NS/Euler system}
\label{mres}

 Let $\mathcal{K}^1 \subset \mathbb{Z}^2$ be a finite  set
of controlled forcing modes; $0 \not\in \mathcal{K}^1$.   Define
the sequence of sets $\mathcal{K}^j \subset \mathbb{Z}^2$
iteratively as:
\begin{eqnarray}\label{setk2}
\mathcal{K}^j  =\mathcal{K}^{j-1} \bigcup \\  \{m+n | \ m,n \in
\mathcal{K}^{j-1} \bigwedge \|m\| \neq \|n\| \bigwedge m \wedge n
\neq 0 \}.  \nonumber
\end{eqnarray}

\begin{thm}{\it (controllability in observed component)}
\label{fulcon}
 Let $\mathcal{K}^1$ be the set of controlled forcing modes, $\mathcal{K}^{obs}$ a finite set of observed
 modes.
 Define iteratively by (\ref{setk2}) the sequence of sets
$\mathcal{K}^j, \ j=2, \ldots ,$ and assume that $\mathcal{K}^M
\supseteq \mathcal{K}^{obs}$ for some $M \geq 1$. Then for any
$T>0$ the 2D NS/Euler system
(\ref{2con1})-(\ref{2cobs})-(\ref{2unobs}) is time-$T$ globally
controllable in the observed component. $\Box$
\end{thm}

%%%This Theorem is proven in Section~\ref{cr}.

\begin{rem}
The present sufficient criterion differs from the one,  we
obtained in \cite{AS43,ASDAN}, by the presence of the 'term'
$\mathcal{K}^{j-1}$ in the right-hand side of the formula
(\ref{setk2}). With the new augmented $\mathcal{K}^j$'s and with
new 'saturation property' (see Definition~\ref{sats})
controllability can be established under weaker hypothesi. $\Box$
\end{rem}

Theorem~\ref{fulcon} characterizes controllability in projection
on a finite-di\-men\-sional {\it coordinate} subspaces. A natural
question is whether the system is controllable in projection on
{\it any} finite-dimensional subspace (a relevance of this
question for regularity of solutions for stochastically forced 2D
NS system has been explained to us by J.C.Mattingly and
E.Pardoux).

\begin{defn}
\label{sats} A finite set $\mathcal{K}^1 \subset \mathbb{Z}^2
\setminus \{0\}$ of forcing modes is called saturating if
$\bigcup_{j=1}^\infty \mathcal{K}^j=\mathbb{Z}^2 \setminus \{0\}$,
where $\mathcal{K}^j$ are defined by (\ref{setk2}). $\Box$
\end{defn}

\begin{thm}{\it (controllability in finite-dimensional
projection)} \label{procon}
 Let $\mathcal{K}^1$ be a saturating  set of controlled forcing
 modes and $\mathcal{L}$ be any finite-dimensional subspace of $H_2$.
  Then for any  $T>0$ the 2D NS/Euler system
(\ref{2con1})-(\ref{2cobs})-(\ref{2unobs}) is time-$T$ solidly
controllable in
 any finite-dimensional  projection. $\Box$
\end{thm}

Another controllability result holds under similar assumptions.

\begin{thm}{\it ($L_2$-approximate controllability)}
\label{apcon}
 Consider the 2D NS/Euler system controlled
by degenerate forcing.  Let $\mathcal{K}^1$ be a saturating  set
of controlled forcing modes.
  Then for any $T>0$ the system
(\ref{2con1})-(\ref{2cobs})-(\ref{2unobs}) is time-$T$
$L_2$-approximately controllable. $\Box$
\end{thm}

\subsection{Saturating sets of forcing modes} \label{satur}

As we see the saturating property is crucial for controllability.
In \cite{HM} a characterization of this property for symmetric
subsets $\mathcal{K}^1 \subset \mathbb{Z}^2$ (such that $k \in
\mathcal{K}^1 \Rightarrow -k \in \mathcal{K}^1$) has been
established.

\begin{prop}[\cite{HM}]\label{minsat}
If a symmetric set  $\mathcal{K}^1 \subset \mathbb{Z}^2$ contains
two vectors which are not collinear and have different lengths,
then $\mathcal{K}^1$ is saturating. $\Box$
\end{prop}

\begin{cor}
The set $\mathcal{K}^1=\{(1,0),(-1,0),(1,1),(-1,-1)\} \subset
\mathbb{Z}^2$ is saturating. $\Box$
\end{cor}

This proves that controllability can be achieved by forcing $4$
modes.

\begin{thm}
Consider 2D NS/Euler system controlled by degenerate forcing. Let
the set of controlled modes contain  vectors $k,\ell,-k,-\ell \in
\mathbb{Z}^2$, where $k \wedge \ell \neq 0, \ \|k\| \neq
\|\ell\|$. Then this system is solidly controllable in any
finite-dimensional projection and is  $L_2$-approximately
controllable. $\Box$
\end{thm}

\section{Relaxation of forcing for 2D NS/Euler system: approximation results
and uniform bounds for trajectories} \label{cont}

In this Section we formulate some results on boundedness
 and continuity of solutions of 2D NS/Euler
system with respect to the forcing. We assume the space of
degenerate forcings to be endowed with a weak topology determined
by so-called relaxation metric. These results are used  in
Section~\ref{cr} for proving controllability in observed
projection. Besides they are interesting for their own sake as  an
example of application of relaxed controls to NS/Euler and other
classes of PDE systems. The proofs are rather technical; they are
to be found in Appendix.

\subsection{Relaxation metric}
\label{rlxmet}

\begin{defn}\label{rxm}
 (see e.g. \cite{G60,Gam}) The relaxation pseudometric
 in the space
\linebreak
 $L^1\left([0,T], \mathbb{R}^d\right)$ is defined by the seminorm
$$ \|u(\cdot)\|_{rx}=\max_{t \in [0,T]}\left\{\left\|\int_0^t u(\tau) d\tau
\right\|_{R^d}\right\}.$$ The relaxation metric is obtained by
 identification of the functions which coincide for almost
all $\tau \in [0,T]. \ \Box$
\end{defn}

The relaxation metric is weaker than the natural metric of
$L^1\left([0,T], \mathbb{R}^d\right)$. The relaxation norms of
fast oscillating functions are small, while their $L_1$-norms can
be large.  For example
$$\|\omega^{1/2} \cos \omega t\|_{rx}=\max_{t \in [0,T]}\left|\int_0^t \omega^{1/2} \cos \omega \tau d\tau
\right| \leq \omega^{-1/2},$$ and $\|\omega^{1/2} \cos \omega
t\|_{rx} \rightarrow 0, \ \mbox{as} \  \omega \rightarrow
+\infty$, while $ \|\omega^{1/2}\cos \omega t\|_{L_1} \rightarrow
+\infty, \ \mbox{as} \  \omega \rightarrow +\infty.$

\begin{lem}
\label{ss:105} Let for integrable functions $\phi_n(\cdot), \
n=1,2, \ldots ,$ their relaxation norms
$\|\phi_n(\cdot)\|_{rx}\stackrel{n \rightarrow
\infty}{\longrightarrow}0$. Let $\{r_\beta(t)| \ \beta \in
\mathcal{B}\}$ be a family of absolutely-continuous functions with
their $W_{1,2}$-norms equibounded:
$$\exists C: \ \|r_\beta(0)\|^2+\int_{\ut}^{\bt}
(\dot{r}_\beta(\tau))^2 d\tau \leq C^2,  \ \forall \beta \in
\mathcal{B}. $$ Then
$\|r_\beta(\cdot)\phi_n(\cdot)\|_{rx}\stackrel{n \rightarrow
\infty}{\longrightarrow}0$, uniformly with respect to $\beta \in
B. \ \Box$
\end{lem}

{\it Proof}.
\begin{eqnarray*}
\left|\int_0^\tau r_\beta(t)\phi_n(t) dt
\right|=\left|r_\beta(\tau)\int_0^\tau \phi_n(t) dt - \int_0^\tau
\dot{r}_\beta(t)\int_0^t\phi_n(\theta) d\theta dt
\right| \leq   \\
  \leq C\left(1+2\sqrt{\tau}\right)\left\|\phi_n(t)\right\|_{rc}.
  \ \Box
\end{eqnarray*}

\subsection{Boundedness of solutions of forced 2D NS/Euler system}
\label{bft}

Consider a set $\mathbf{F}$ of degenerate forcings $v(t,x)=\sum_{k
\in \mathcal{K}^1}v_k(t)e^{ik \cdot x}$;  $\sharp
\mathcal{K}^1=d$. We identify these forcings with vector-functions
$v(t)=\left(v_k(t)\right) \in L_\infty([0,T];\mathbb{R}^d)$.
Forced 2D NS/Euler system is treated as an evolution equation in
$H_s, \ s \geq 2$. An example of boundedness result, we are
interested in, would be the following

\begin{lem}
\label{forc} Assume the set $\mathbf{F}$ of degenerate forcings to
be bounded in the relaxation metric. Fix the time interval $[0,T]$
and the initial condition $w(0)=w_0 \in H_s, \ s \geq 2$ for the
2D NS/Euler system. Then the tra\-jec\-to\-ries $w_t$ of the
system (\ref{nsw}) forced by $v(t,x) \in \mathbf{F}$  are
equibounded in $H_0$ norm:
$$\exists b: \ \mbox{vrai}\sup_{t \in [0,T]}
\|w_t\|_{0} \leq b. \ \Box $$
\end{lem}

This result is not covered by classical results on boundedness of
solutions of 2D NS/Euler system because the set of forcings can be
bounded in the relaxation metric while being unbounded in
$L_\infty$ and $L_2$ metric.

We will derive the previous result from a stronger assertion
(Theorem~\ref{bound}). To formulate the assertion  consider the
primitives $V(\cdot)=\int_0^\cdot v(\tau) d \tau$ of $v(\cdot) \in
\mathbf{F}$. By assumptions of the Lemma~\ref{forc} $V(\cdot)$ are
equibounded in the metric of $C^0([0,T],\mathbb{R}^d)$.

Denote the trigonometric polynomial $\sum_{k \in
\mathcal{K}^1}V_k(t)e^{ik \cdot x}$ by $V_t(x)$. The forced
(controlled) 2D Euler system can be written as
\begin{equation}\label{dvte}
  \partial w_t/ \partial t= (u_t \cdot \nabla)w_t+
\nu\Delta w_t+  \partial V_t / \partial t .
\end{equation}

Put $y_t=w_t-V_t$.
 The equation~(\ref{dvte}) can be rewritten as:
\begin{equation}\label{dyte}
\partial y_t/ \partial t = (u_t \cdot
\nabla)(y_t+V_t)+\nu \Delta y_t+\nu \Delta V_t  .
\end{equation}

Recall that $u_t$  is the divergence-free solution of the equation
$ \nabla^\perp \cdot u_t=w_t=y_t+V_t$. It can be represented as a
sum ${\mathcal Y}_t+{\mathcal V}_t$, where ${\mathcal
V}_t,{\mathcal Y}_t$ are the divergence-free solutions of the
equations $ \nabla^\perp \cdot {\mathcal V}_t=V_t,\ \nabla^\perp
\cdot {\mathcal Y}_t=y_t,$  with periodic boundary conditions.

Hence the equation (\ref{dyte}) allows for the representation
\begin{eqnarray}\label{yvv}
 \partial y_t/\partial t=\left(\left({\mathcal Y}_t+{\mathcal V}_t\right)
    \cdot \nabla \right) \left(y_t+V_t\right) \\ =\left({\mathcal Y}_t
    \cdot \nabla \right) y_t+\left({\mathcal V}_t
    \cdot \nabla \right) y_t+\left({\mathcal Y}_t
    \cdot \nabla \right) V_t+\nu \Delta y_t+\nu \Delta V_t +\left({\mathcal V}_t
    \cdot \nabla \right) V_t.  \nonumber
\end{eqnarray}

This equation can be seen as 2D NS/Euler equation forced by
$y$-linear forcing term $\left({\mathcal V}_t
    \cdot \nabla \right) y_t+\left({\mathcal Y}_t
    \cdot \nabla \right) V_t $ together with $y$-independent forcing term
 $\nu \Delta V_t +\left({\mathcal V}_t
    \cdot \nabla \right) V_t$.

Consider instead of (\ref{yvv})  a more general equation
\begin{equation}\label{ef}
\partial y_t/\partial t=\left({\mathcal Y}_t
    \cdot \nabla \right) y_t+\left({\mathcal V}^1_t
    \cdot \nabla \right) y_t+\left({\mathcal Y}_t
    \cdot \nabla \right) V^2_t+ \nu \Delta y_t+V^0_t,
\end{equation}
where all the forcing terms $V^0_t,V^1_t,V^2_t$ are now decoupled
and  ${\mathcal V}^1_t$ is the divergence-free solution of the
equation $\nabla^\perp \cdot {\mathcal V}^1_t=V^1_t$.

Consider the set
 $\mathbf{F}_B$  of triples $\{(V^0_t,V^1_t,V^2_t)\}$ satisfying the
 condition:
\begin{equation}\label{bv012}
 \sup_{t \in [0,T]}\max\{\|V^0_t\|, \| V^1_t\|, \|V^2_t\|\} \leq
B, \ B>0.
\end{equation}
As far as $V^i_t$ are trigonometric polynomials (of fixed order)
in $x$ the norms $\|\cdot\|_{H_s}$ are equivalent for all $s$, so
one just uses the notation $\| \cdot \|$.

All the results of this Section will be proven for the equation
(\ref{ef}) and for the forcings from the set $\mathbf{F}_B$
defined by (\ref{bv012}).

\begin{rem}
The equation (\ref{ef}) is "nonclasically forced" 2D NS/Euler
equation.  Remarks on existence and uniqueness results for its
solutions can be found in Appendix. $\Box$
\end{rem}

\begin{thm}
\label{bound} Let  $\mathbf{F}_B=\{(V^0_t,V^1_t,V^2_t)\}$ be the
set defined by (\ref{bv012}). Fix the time interval $[0,T]$ and
the initial condition $y(0)=y_0 \in H_s, \ s \geq 2$ for the
system (\ref{ef}) forced by the elements of $\mathbf{F}_B$. Then
$\exists b>0$ such that for all $(V^0_t,V^1_t,V^2_t) \in
\mathbf{F}_B$ and for the corresponding trajectories $y_t$ of the
equation (\ref{ef}) there holds:

\begin{equation}\label{eqblin}
\mbox{i)}  \ vrai\sup_{t \in [0,T]} \|y_t\|_{L_\infty} \leq b;
\end{equation}
\begin{equation}\label{eqbh1}
\mbox{ii)}  \ vrai\sup_{t \in [0,T]} \|y_t\|_{H_2} \leq b;
\end{equation}
\begin{equation}\label{intdy1}
 \mbox{iii)}\  \int_0^T
\left\|\frac{\partial}{\partial t}y_t\right\|^2_1 dt \leq b. \
\Box
\end{equation}
\end{thm}

\begin{rem}
Obviously the conclusion of the Lemma~\ref{forc} can be derived
from this theorem.
\end{rem}

The proof of the Theorem~\ref{bound} is to be found in the
Appendix.

%%%%%%%%%%%%%%%%%%%%%%%%%%%%%%%%%%%%%%%%%%%%%%%%%%%%%%%%%%%%%%%%%%%%
\subsection{Continuous dependence of tra\-jec\-to\-ries on relaxed forcings}
\label{cft}

In this subsection we establish continuous dependence of
trajectories of the equation (\ref{ef}) on the forcing terms
$V^1_t,V^2_t,V^0_t$, as these latter vary continuously in the
relaxation metric.

\begin{thm}\label{corr} Consider the set
$\mathbf{F}_B=\{(V^0_t,V^1_t,V^2_t)\}$ defined by (\ref{bv012}).
Fix the time interval $[0,T]$ and the initial condition $y(0)=y_0
\in H_s, \ s \geq 2$ for the system (\ref{ef}) forced by the
elements of $\mathbf{F}_B$. Endow $\mathbf{F}_B$ with the
relaxation metric and endow the space of trajectories of the 2D
NS/Euler equation with $L_\infty((0,T);H_0)$-metric. Then the
restriction of  the forcing/trajectory map onto $\mathbf{F}_B$ is
uniformly continuous. $\Box$
\end{thm}

The proof of the Theorem~\ref{corr} is to be found in Appendix.

\section{Proof of controllability in observed component for 2D
NS/Euler system}
\label{cr}

We prove first the result on solid controllability in observed
component (Theorem~\ref{fulcon}). The construction introduced in
proof is used to establish controllability  in finite-dimensional
projection and  $L_2$-approximate controllability.

First we slightly particularize the assertion of the
Theorem~\ref{fulcon}.
\begin{thm}
\label{fam} Let $\mathcal{K}^1$ be a set of controlled forcing
modes. Define according to (\ref{setk2}) the sequence of sets
$\mathcal{K}^j, \ j=2, \ldots ,$
 and assume that, for some $M$,
 $\mathcal{K}^M \supset \mathcal{K}^{obs}$.

Then for all sufficiently small $T>0$ the 2D NS/Euler system is
solidly controllable in projection on the observed component.
Besides one can choose the corresponding family of controls
$v(\cdot ,b)$ (cf. Definition~\ref{socoproj}) which is
parameterized continuously in $L_1$-metric by a compact subset
$B_R$ of a finite-dimensional linear space and is uniformly (with
respect to $t,b$) bounded: $\forall t,b: \ \|v(t,b)\| \leq A(T,R).
\ \Box$.
\end{thm}

The only additional restriction in the claim of the latter result
is smallness of time.  To deal with large $T$ we can apply zero
control on the interval $[0,T-\theta]$ with $\theta$ small and
then apply the result of the Theorem~\ref{fam}.

\subsection{Sketch of the proof}

By assumption the set $\mathcal{K}^{obs}$ of observed modes is
contained in some $\mathcal{K}^M, \ M \geq 1,$ from the sequence
defined by (\ref{setk2}). We will proceed by induction on $M$.

If $M=1$ then $\mathcal{K}^1 \supset \mathcal{K}^{obs}$, i.e. all
the  equations for the observed modes contain controls. Then it is
easy to establish small time controllability in observed
component, given the fact that there are no a priori bounds on
controls (this is done in Subsection~\ref{stepm1}).

Let  $M>1, \ \mathcal{K}^M \supset \mathcal{K}^{obs} \supset
\mathcal{K}^1$. We start acting as if independent control
parameters enter all the equations indexed by $k \in
\mathcal{K}^M$. Then we are under previous assumption and hence
can construct a needed family of controls. Though the control
parameters indexed by $k \in \mathcal{K}^M \setminus
\mathcal{K}^1$ are fictitious and our next step would be
approximating the actuation of (some of) these fictitious controls
by actuation of controls of {\it smaller} dimension. Now we employ
the controls, which only enter the equations indexed by $k \in
\mathcal{K}^{M-1} \subset \mathcal{K}^M$. The possibility of such
approximation for 2D NS/Euler system (provided that the relation
between $\mathcal{K}^{M-1}$ and $\mathcal{K}^M $ is established by
(\ref{setk2})) is the main element of our construction.

If $M-1>1$, then  the approximating  controls are also fictitious,
but we can repeat the reasoning in order to arrive  after $M-1$
steps to {\it true} controls indexed by $\mathcal{K}^1$.

We can look at the process the other way around. Starting with a
(specially chosen) family of degenerate controls in low modes,
indexed by $\mathcal{K}^1$  we transfer their actuation to the
higher modes via the {\it nonlinear term} of 2D NS/Euler system.

\subsection{Proof of the Theorem~\ref{fam}: first induction
step} \label{stepm1}

 The first induction step ($M=1$) follows from the
following Lemma.

\begin{lem}
\label{m=1}
 Let $M=1$, and $\mathcal{K}^1 = \mathcal{K}^{obs}$. The 2D NS/Euler system is split in the subsystems (\ref{2con1}) and
 (\ref{2unobs}),
 which can be written in a concise form as
\begin{eqnarray}\label{q1Q}
  d q^1/dt=f_1(q^1,Q) + v, \
d Q/d t= F(q^1,Q)\\
q^1(0)=q^1_0, \ Q(0)=Q_0, \nonumber
\end{eqnarray}
$\dim q^1=N$.  Then for sufficiently small $\tau>0$:
   there exists a family of controls  $v(t;b)$ which satisfies
   the conclusion of the Theorem~\ref{fam}. $\Box$
\end{lem}

{\em Proof.} Without lack of generality we may assume the initial
condition for the observed component  to be
$q^1(0)=0_{R^{\kappa_1}}$. We do not diminish generality either by
assuming $\mathcal{K}^1 = \mathcal{K}^{obs}$ instead of
$\mathcal{K}^1 \supseteq \mathcal{K}^{obs}$. Recall that
$\Pi_1:(q^1,Q) \rightarrow q^1$.

Define for $y \in \mathbb{R}^N $,  $\|y\|_{l_1}=\sum_{j=1}^N
|y_i|$. Let $\mathcal{C}_R=\{y \in \mathbb{R}^N| \ \|y\|_{l_1}
\leq R\}$.

Fix $\gamma >1$. Take the interval $[0,\tau]$; the value of {\em
small} $\tau >0$ will be specified later on.  For each $b \in
\gamma \mathcal{C}_R$ take $v(t;p,\tau) =\tau^{-1}p$ - a constant
control. Obviously $\gamma \mathcal{C}_R \supset \mathcal{C}_R$
and $\int_0^{\tau} v(t;p,\tau)dt=p$. For fixed $\tau >0 $ the map
$p \mapsto v(t;p,\tau)$ is continuous in $L_1$-metric.

We claim that $\exists \tau_0>0$ such that for $\tau \in (0,
\tau_0)$ the family of controls $v(t;p,\tau), \ p \in \gamma
\mathcal{C}_R$ satisfies the conclusion of the Lemma, so one may
take $b=p, \ B_R=\gamma \mathcal{C}_R$.

Denote for fixed $\tau > 0$ the map
$$p \mapsto v(\cdot ;p,\tau) \mapsto (\Pi_1 \circ
\mathcal{F/O}_{\tau})\left(v(\cdot;p,\tau)\right)$$ by
$\Phi(p;\tau)$.  Recall that $\Pi_1 \circ \mathcal{F/O}_{\tau}$ is
the end-point component map (cf. Definition~\ref{itomap}). The map
$p \mapsto v(\cdot;p,\tau)$ is continuous in $L_1$-metric of
controls and hence also in the relaxation metric. Therefore by
Theorem~\ref{corr} the map $p \mapsto \Phi(p;\tau)$ is continuous.

Restrict the equations (\ref{q1Q}) to the interval $[0,\tau]$ and
proceed with time substitution $t=\tau \xi, \ \xi \in [0,1]$. The
equations  take form:
\begin{equation}\label{resc}
 d q^1/d\xi=\tau f_1(q^1,Q) + p, \
d Q/d \xi = \tau F(q^1,Q), \ \xi \in [0,1].
\end{equation}

For $\tau=0$ the 'limit system' of (\ref{resc}) is
\begin{equation}\label{0resc}
   d q^1_0/d\xi= p, \
d Q_0/d \xi = 0, \ \xi \in [0,1].
\end{equation}

The end-point component map $p \mapsto q^1_0(1)$ for the limit
system is the identity.

From classical results on boundedness of solutions of 2D NS/Euler
system we conclude that
 $q^1$-components of the solutions of the systems (\ref{resc}) and (\ref{0resc})
 (with the same initial
condition)  deviate
 by a quantity $\leq C\tau$,
where the constant $C$ can be chosen independent of $p,\tau$ for
sufficiently small $\tau >0$. Then $\|\Phi(\cdot;\tau)-Id\| \leq C
\tau$. By degree theory argument there exists $\tau_0$ such that
$\forall \tau \in (0,\tau_0)$ the image of $p \mapsto
\Phi(p;\tau)$ covers $\mathcal{C}_R$ solidly.

To complete the proof note that $\|v(t;p,\tau)\|$ are uniformly
bounded by $\gamma R \tau^{-1}. \ \Box$

In what follows we will need a modification of the previous Lemma.

\begin{lem}\label{m1par}
Consider the system (\ref{q1Q}) and impose the boundary conditions
$$q^1(0)=\phi(p), q^1(\tau)=\psi(p), \ p \in P, \ P \ \mbox{-
compact}, \ \phi , \psi  \ \mbox{- continuous},$$ on its
$q^1$-component.

Then for all sufficiently small $\tau>0$, there exists a family of
controls $v(t;p,\tau)$ defined on $[0,\tau]$, such that the
corresponding trajectories, which meet the initial condition, meet
the end-point condition approximately:
$$\|q^1(\tau;p)-\psi(p)\| \leq C\tau.$$
 Besides $\|Q(t)-Q_0\|_0 \leq \gamma C\tau, \ \forall t \in [0,\tau]$. Here $C$ can be chosen independent on $p , \tau . \ \Box$
\end{lem}

The proof is similar to the previous one. One can choose the
family of controls $v(t,p)=\tau^{-1}(\psi(p)-\phi(p)), \ t \in
[0,\tau]$.

\subsection{Generic induction
step: solid controllability by extended controls} \label{stepgen}

Let us proceed further with the induction. Assume that the
statement of the Theorem~\ref{fam} has been proven for all $M \leq
(N-1)$; we are going to prove it for $M=N$.

  Consider now the  system
(\ref{2con1})-(\ref{2cobs})-(\ref{2unobs}) with the {\it extended}
set $\mathcal{K}^1_e=\mathcal{K}^2$ of controlled forcing modes;
$\mathcal{K}^1_e \supseteq \mathcal{K}^1$.

Obviously this new system satisfies the conditions of the
Theorem~\ref{fam}; indeed
$$\mathcal{K}^1_e=\mathcal{K}^2 \Rightarrow
\mathcal{K}^j_e=\mathcal{K}^{j+1}, \ j \geq 1,$$ for the sets
$\mathcal{K}^j_e, \mathcal{K}^j$, defined by (\ref{setk2}). Hence
$\mathcal{K}_e^{M-1}=\mathcal{K}^M \supseteq \mathcal{K}^{obs}$.

By  induction hypothesis the system with extended controls is
solidly controllable in observed projection: there exists a
continuous in $L_1$-metric family of extended controls $v(t;b)$
which satisfies the conclusion of the Theorem~\ref{fam}.

This family of controls is uniformly bounded; assume that $
\|v(t;b)\|_{l_1} \leq A, \ \forall b \in B, \ \forall t \in
    [0,T]$.
The values of $v(t;b)$ belong to $\mathbb{R}^{\kappa_2}$, where
$\kappa_2=\# \mathcal{K}^2=\mathcal{K}^1_e$.

Evidently these extended controls are unavailable for the original
problem. We are going  to approximate their action by the action
controls from a more restricted set.

To this end let us first take the vectors $e_1, \ldots , e_{k_2}$
from the standard basis in $\mathbb{R}^{k_2}$ together with their
opposites $-e_1, \ldots , -e_{k_2}$.   Multiply each of these
vectors by $A$ and denote the set of these $2\kappa_2$ vectors  by
$E^A_2$. The convex hull $\mbox{conv} E^A_2$ of $E^A_2$ contains
all the values of
 $v(t;b)$.

First we  will approximate the family of functions $v(t;b)$ which
take their values in $\mbox{conv} E^A_2$ by $E^A_2$-valued
functions. Such a possibility is a central result of relaxation
theory.

\begin{defn}
 Define $\delta$-pseudometric $\rho_\delta$ in the space
 $L^\infty\left([0,T],\mathbb{R}^\kappa \right)$ of measurable functions
 as:
 $$\rho_\delta\left(u^1(\cdot),u^2(\cdot)\right)=\mbox{meas}\{t \in [0,T] \left|  \right. u^1(t) \neq
 u^2(t)\}.$$ Identifying the functions, which coincide beyond
 a set of zero measure, we arrive to $\delta$-metric. $\Box$
\end{defn}

\begin{rem}
The $\delta$-metric is a restriction of strong metric of relaxed
controls (see \cite{Gam}) onto the  set of ordinary (=nonrelaxed)
controls. $\Box$
\end{rem}

We will apply  R.V.Gamkrelidze Approximation Lemma (see
\cite[Ch.3]{Gam},\cite[p.119]{G60}). According to it given a
$\delta$-continuous family of $\mbox{conv } E_2$-valued functions
and $\eps
>0$ one can construct  a $\delta$-continuous family
 of $E_2$-valued functions which $\eps$-approximates the
family $\{v(t;b) | b \in B\}$ in the relaxation metric uniformly
with respect to $b \in B$. Moreover the functions of the family
can be chosen piecewise-constant and the number $L$ of the
intervals of constancy can be chosen the same for all $b \in B$.
Actually the Approximation Lemma in \cite[Ch.3]{Gam} regards
relaxed controls (Young measures). Applying  it to nonrelaxed
controls $v(\cdot;b)$ we arrive to the following result.

\begin{prop}(cf. Approximation Lemma; \cite[Ch.3]{Gam}).
\label{alem}  Let $B$ be a compact and $\{v(t;b) | b \in B\}$ be a
family of $(\mbox{conv }E^A_2)$-valued functions, which depends on
$b \in B$ continuously in $L_1$ metric. Then for each $\eps >0$
one can construct a $\delta$-continuous (and hence
$L_1$-continuous) equibounded family $\{z(t;b)| \ b \in B\}$ of
$E^A_2$-valued functions which $\eps$-approximates the family
$\{v(t;b) | b \in B\}$ in the relaxation metric uniformly with
respect to $b \in B$. Moreover the functions $z(t;b)$ can be
chosen piecewise-constant and the number $L$ of the intervals of
constancy can be chosen the same for all $b \in B$. The intervals
of constancy of these controls vary continuously with  $b \in B \
\Box$
\end{prop}

We omit the proof, which is a slight variation of the proof in
\cite[Ch.3]{Gam}.

Applying this result to our case we construct a
$\delta$-continuous family of $E^A_2$-valued functions $\{z(t;b)|
\ b \in B\}$ which approximates the family $\{v(t;b) | \ b \in
B\}$ uniformly in the relaxation metric. According to the
Theorem~\ref{corr} the end-point map $\mathcal{F/O}_T$ is
continuous in the relaxation metric. Therefore we conclude with
the following result.

\begin{prop}
\label{decon}  There exist a number $L$ and an $L_1$-continuous
family of piece\-wise-constant $E^A_2$-valued controls $\{z(t;b)|
\ b \in B\}$  (with at most $L$ intervals of constancy) such that
the reduced system is solidly controllable by means of this
family. $\Box$
\end{prop}

\subsection{Generic induction
step: solid controllability of the original system}
\label{stepfin}

Let us compare the original system
(\ref{2con1})-(\ref{2cobs})-(\ref{2unobs}) with the system driven
by the $E^A_2$-valued controls $\{z(t;b)| \ b \in B\}$ constructed
in the Proposition~\ref{decon}.

In both systems the equations for the coordinates $q_k$, indexed
by $k \in {\mathcal K}^1$ coincide:
\begin{equation}\label{eq:5:3}
\dot{q}_k=\sum_{m+n=k}(m \wedge n)|m|^{-2}q_m q_n-\nu|k|^2q_k+v_k,
\ k \in {\mathcal K}^1.
\end{equation}
We collect these coordinates into the vector denoted by $q^1$.

In the original system the equations for the variables $q_k, \ k
\in \left({\mathcal K}^2 \setminus {\mathcal K}^1\right)$ are
'uncontrolled':
\begin{equation}\label{k2-1}
 \dot{q}_k=\sum_{m+n=k}(m \wedge n)|m|^{-2}q_m
q_n - \nu|k|^2q_k, \ k \in \left({\mathcal K}^2 \setminus
{\mathcal K}^1\right).
\end{equation}
They differ from the corresponding equations of the system with
extended controls, which are:
\begin{equation}\label{eq:k2}
  \dot{q}_k=\sum_{m+n=k}(m \wedge n)|m|^{-2}q_m q_n- \nu|k|^2q_k+z_k, \ k \in \left(\mathcal{K}^2 \setminus {\mathcal
K}^1\right).
 \end{equation}
 We collect $q_k, \ k \in \left(\mathcal{K}^2 \setminus
{\mathcal K}^1\right)$ into the vector denoted by $q^2$ and denote
$q=(q^1,q^2)$.

Finally the equation for the infinite-dimensional component $Q_t$,
which collects the higher modes $e^{i k\cdot x}, \ k \not\in
{\mathcal K}^2,$ does not contain controls and is the same in both
systems. It suffices for our goals to write this equation in  a
concise form as:
\begin{equation}\label{2perp}
\dot{Q}=h(q,Q).
\end{equation}

According to the Proposition~\ref{decon} we manage to control our
system solidly by means of extended $E^A_2$-valued
piecewise-constant controls $z(t;b)$. By Proposition~\ref{alem}
the intervals of constancy vary continuously with $b \in B$. Our
task now is to design a family of "small-dimensional" controls
$x(t;b)$ for the equations
(\ref{eq:5:3})-(\ref{k2-1})-(\ref{2perp}), such that the maps
$$b
\mapsto z(\cdot ;b) \mapsto (\Pi^{obs} \circ
\mathcal{F/O}_{\tau})\left(z(\cdot;b)\right) \ \mbox{and} \ b
\mapsto x(\cdot ;b) \mapsto (\Pi^{obs} \circ
\mathcal{F/O}_{\tau})\left(x(\cdot;b)\right)$$ are $C^0$-close.

If on some interval of constancy $T \in [\ut,\bt]$ the value of
$z(t;b)$ equals $\pm A e_k$ with $k \in {\mathcal K}^1$, then we
just take the control $x(t;b)$ in (\ref{eq:5:3}) coinciding with
$z(\cdot ;b)$ on this interval.

The real problem arises when on some interval of constancy
$z(t;b)$ takes value  $\pm A e_{\bar{k}}$ with $\bar{k} \in
\left(\mathcal{K}^2 \setminus \mathcal{K}^1 \right)$. There are no
controls available in the corresponding equation (\ref{k2-1}) for
$q_{\bar{k}}$ and we will "affect" the evolution of $q_{\bar{k}}$
via the variables $q_m, \ m \in \mathcal{K}^1$ which enter this
equation.

More exactly the construction of the controls $x(t;b,\omega)$ on
the intervals of constancy of $z(t;b)$ goes as follows:

 1) on an interval $[\ut,\bt]$ of the first kind, where $z(t;b)=e_k$ with $k
\in {\mathcal K}^1$ we take $x(t;b,\omega) =z(t;b)$;

2) on an interval $[\ut,\bt]$ of the second kind, where
$z(t;b)=A(b)e_k$ with $k \in {\mathcal K}^2 \setminus {\mathcal
K}^1$, pick a pair $m,n \in \mathcal{K}^1$ such that $$m \wedge n
\neq 0, \ |m| \neq |n|, \ m+n=k,$$ such  pair exists by the
definition of ${\mathcal K}^2$ (see (\ref{setk2})); choose
$A_m(b),A_n(b)\in \mathbb{R}$ satisfying
\begin{equation}\label{Amn}
|A_m(b)|=|A_n(b)|   \bigwedge A_m(b) A_n(b)(m \wedge
n)(|m|^{-2}-|n|^{-2})=2A(b).
\end{equation}
and take
$$x_m(t;b,\omega)=A_m (b)\omega \cos \omega t, \ x_n(t;b,\omega)=A_n (b)\omega \cos \omega t,$$
choosing other components  $x_j(t;b,\omega)$ equal to $0$.

It is easy to see that the  primitives $X(t;b,\omega)=\int_0^t
x(s;b,\omega)ds$ are bounded by a constant which can be chosen
independent of $b$ and $\omega$. Besides $X(T;b,\omega)$ varies
continuously with $b$ (for fixed $\omega$).

%%%%%%%%%%%%%%%%%%%%%%%%%%%%%%%%%%%%%%%%%%%%%%%%%%%%%

Consider two   trajectories $w^{b,\omega}_t, \bar{w}_t, \ t \in
[0,T]$, which are  driven by the controls $x(t;b,\omega)$ and
$z(t;b)$ correspondingly. We will prove that $w^{b,\omega}_t$ and
$\bar{w}_t$ match asymptotically (as $\omega \rightarrow \infty$)
in all the components but $q^1$. Let $\Pi_1$ be the projection
onto the space of modes $\{e^{im\cdot x}| \  m \in {\mathcal
K}^1\}$, while $\Pi_1^\perp$ be the projection onto its orthogonal
complement.

%%%%%%%%%%%%%%%%%%%%%%%%%%%%%%%%%%%%%%%%%%%%%%%%%%%
\begin{prop}
\label{10.5} i) The trajectories $w^{b,\omega}_t$ are equibounded:
$$\exists C: \ \|w^{b,\omega}_t\|_0 \leq C, \ \forall t \in [0,T],
\ b \in B, \ \omega >0;$$

ii) For fixed $\omega >0$ the dependence $b \mapsto
w^{b,\omega}_t$ on $b$ is continuous in $C^0[0,T]$-metric of the
controls;

iii)  For any $\eps >0$ there exists $\delta >0$ and $\omega_0$
such that if $\omega
> \omega_0$ and
$\|w^{b,\omega}|_{t=0}-\bar{w}|_{t=0}\|_0 \leq \delta$, then
$\forall t \in[0,T]: \
\|\Pi_1^\perp\left(w^{b,\omega}_{t}-\bar{w}_{t}\right)\|_0 \leq
\eps. \ \Box$
\end{prop}

Assuming the claim of this Proposition (which is proven in the
Appendix) to hold true let us complete the induction.

By assumption the system
(\ref{eq:5:3})-(\ref{eq:k2})-(\ref{2perp}) is solidly controllable
in observed component by means of the family of extended controls
$z(t;b)$, i.e. the map $b \mapsto \left(\Pi^{obs} \circ
\mathcal{F/T}_T\right)(z(t;b))$ covers solidly the cube
$\mathcal{C}_R$ in $\Pi^{obs}(H_2)$. According  to the
Proposition~{\ref{10.5} all the components, but $q^1$, of the
trajectories driven by $x(t;b,\omega_0)$ match up to arbitrarily
small $\eps$, provided $\omega$  is sufficiently large. For fixed
$\omega$ the controls $x(t;b,\omega)$ depend continuously (in
$L_1$-metric) on $b \in B$.

By the degree theory argument we may conclude  that for large
$\omega$ the map
$$b
\mapsto \left(\Pi_1^\perp \circ \Pi^{obs} \circ
\mathcal{F/T}_T\right)(x(t;b,\omega))$$
 covers solidly the set
$\Pi_1^\perp(\mathcal{C}_R)$.

Still the map  $b \mapsto (\Pi_1 \circ
\mathcal{F/T}_T)(x(t;b,\omega))$ does not necessarily match with
$b \mapsto (\Pi_1 \circ \mathcal{F/T}_T)(z(t;b))$. We have to
settle the $q^1$-component.

Considering $q^1(t;\omega, b)$ and evaluating  it at $T$ we
observe that according to the Proposition~\ref{10.5} the values
$q^1(T;\omega, b)$ are equibounded for all $b, \omega$ and for
fixed $\omega >0$ the dependence $b \rightarrow q^1(T;\omega, b)$
is continuous. Put
$$\phi(b)=q^1(T;\omega, b), \psi(b)=\left(\Pi_1 \circ
\mathcal{F/T}_T\right)(z(t;b)).$$
  We can apply
Lemma~\ref{m1par} for constructing  controls $x(t;b,\omega)$
defined on an arbitrarily small interval $[T, T+\tau]$ such that
\begin{eqnarray*}
\|q^1(T+\tau;b)- \psi(b)\|=O(\tau), \\  \|(\Pi_1^\perp \circ
\mathcal{F/T}_{T+\tau})(x(t;\omega))-(\Pi_1^\perp \circ
\mathcal{F/T}_{T})(z(t;b))\| =O(\tau), \ \mbox{as} \ \tau
\rightarrow 0.
\end{eqnarray*}

 Then choosing $\tau >0$ sufficiently small  we prove
that the maps
$$b \mapsto (\Pi^{obs} \circ
\mathcal{F/T}_T)(z(t;b)) \ \mbox{and}\  b \mapsto (\Pi^{obs} \circ
\mathcal{F/T}_{T+\tau})(x(t;b,\omega))$$  are close in
$C^0$-metric and therefore by the degree theory argument the last
map covers solidly the cube $\mathcal{C}_R$. It means that the
system is time-$(T+\tau)$ solidly controllable. $\Box$

\subsection{Proof of the Proposition~\ref{10.5}}
We proceed by induction on a uniformly (with respect to $b \in B$)
bounded number of the intervals of constancy of controls
$z(\cdot;b)$. Since on the intervals of the first kind the
controls $z(\cdot;b)$ and $x(\cdot;b)$ coincide it suffices to
consider one interval $[\ut,\bt]$ of the second kind. We may think
that $[\ut,\bt]=[0,T]$.

Introduce the primitive $V^\omega_t(b)=\int_0^t
x(\tau;\omega,b)d\tau$  of the control $x(\cdot;b,\omega)$. The 2D
NS/Euler system can be written as
$$\partial w^\omega_t/ \partial t= (u^\omega_t \cdot
\nabla)w^\omega_t+
  \nu \Delta w^\omega_t + \partial V^\omega_t / \partial t .$$

Introduce  $y^\omega_t=w^\omega_t-V^\omega_t$. Notice that
$y^\omega_t$ and $w^\omega_t$ differ only in ${\mathcal
K}^1$-indexed modes, i.e. $\Pi_1^\perp y^\omega_t=\Pi_1^\perp
w^\omega_t$.
 The  equation for $y^\omega_t$ is:
\begin{equation}\label{yom}
\partial y^\omega_t/ \partial t = (u^\omega_t \cdot
\nabla)(y^\omega_t+V^\omega_t)+
  \nu \Delta \left(y^\omega_t+V^\omega_t\right).
\end{equation}

The function $u^\omega_t$  can be represented as a sum ${\mathcal
Y}^\omega_t+{\mathcal V}^\omega_t$, where ${\mathcal
V}^\omega_t,{\mathcal Y}^\omega_t$ are the divergence-free
solutions of the equations: $ \nabla^\perp \cdot {\mathcal
V}^\omega_t=V^\omega_t,\ \nabla^\perp \cdot {\mathcal
Y}^\omega_t=y^\omega_t,$  under periodic boundary conditions.

Hence the equation (\ref{yom}) allows for the representation
\begin{eqnarray*}
 \partial y^\omega_t/\partial t=\left(\left({\mathcal Y}^\omega_t+{\mathcal V}^\omega_t\right)
    \cdot \nabla \right) \left(y^\omega_t+V^\omega_t\right)+\nu \Delta \left(y^\omega_t+V^\omega_t\right)= \left({\mathcal Y}^\omega_t
    \cdot \nabla \right) y^\omega_t+  \\   +\left({\mathcal V}^\omega_t
    \cdot \nabla \right) y^\omega_t+\left({\mathcal Y}^\omega_t
    \cdot \nabla \right) V^\omega_t+\nu \Delta y^\omega_t+ \nu \Delta V^\omega_t+ \left({\mathcal V}^\omega_t
    \cdot \nabla \right) V^\omega_t.
\end{eqnarray*}

 Denote $e^{i\ell \cdot x}$ by $e_\ell$; then $V^\omega_t=
(A_me_m+A_m e_n)\sin \omega t$. Obviously ${\mathcal
V}^\omega_t=(A_m{\mathcal V}_m+A_n{\mathcal V}_n)\sin \omega t$.

On an interval of the second kind the equation for $\bar{w}_t$,
driven by the constant control $z(t;b)$, is:
$$\partial_t \bar{w}_t =\left(\bar{\mathcal W}_t \cdot \nabla
\right) \bar{w}_t+\nu \Delta \bar{w}_t
    +Ae_{m+n}, \ m,n \in \mathcal{K}^1,$$
where $\bar{\mathcal W}_t$ is the divergence-free solution of the
equation $\nabla^\perp \cdot \bar{{\mathcal W}}_t= \bar{w}_t$
 under periodic boundary conditions.

Introducing the notation $\eta^\omega_t=y^\omega_t-\bar{w}_t$, we
obtain for $\eta^\omega_t$ the equations:
$$
\partial_t \eta^\omega_t =\left(\left({\mathcal W}^\omega_t+{\mathcal V}^\omega_t\right)
    \cdot \nabla \right) y^\omega_t-\left(\bar{\mathcal W}_t \cdot \nabla \right)
    \bar{w}_t+ \nu\Delta \eta^\omega_t+\left(\left({\mathcal V}^\omega_t \cdot \nabla
    \right)V^\omega_t-Ae_{m+n}\right).
$$
Subtracting and adding $\left(\left({\mathcal
W}^\omega_t+{\mathcal V}^\omega_t\right) \cdot \nabla \right)
\bar{w}_t$     to the right-hand side of the latter equation we
transform it into
\begin{eqnarray}\label{eom2}
\partial_t \eta^\omega_t =\left({\mathcal H}^\omega_t
    \cdot \nabla \right) \eta^\omega_t+
\left(\left(\bar{\mathcal W}_t+  {\mathcal V}^\omega_t \right)
\cdot \nabla \right) \eta^\omega_t+
    \left(\left({\mathcal H}^\omega_t + {\mathcal V}^\omega_t
    \right)
    \cdot \nabla \right) \bar{w}_t+  \nonumber \\+ \nu\Delta \eta^\omega_t    +\left(\left({\mathcal V}^\omega_t \cdot \nabla
    \right)V^\omega_t-Ae_{m+n}\right),
\end{eqnarray}
where ${\mathcal H}^\omega_t$ is the divergence-free solution of
the equation $\nabla^\perp \cdot {\mathcal
H}^\omega_t=\eta^\omega_t$.

By construction $V^\omega_t, {\mathcal
V}^\omega_t,\left(\left({\mathcal V}^\omega_t \cdot \nabla
\right)V^\omega_t-Ae_{m+n}\right)$ converge to $0$ in relaxation
metric, as $\omega \rightarrow +\infty$. By the continuity result
(Theorem~\ref{corr})  trajectories of (\ref{eom2}) converge in
$C^0$  to the trajectories of the equation
\begin{equation}\label{limeq}
\partial_t \eta_t =\left(\left({\mathcal
H}_t+\bar{\mathcal W}_t\right)
    \cdot \nabla \right) \eta_t+\left({\mathcal H}_t
    \cdot \nabla \right) \bar{w}_t+\nu\Delta \eta_t,
\end{equation}
as $\omega \rightarrow +\infty$,

To  estimate  the evolution of $\|\eta_t\|_0$ by virtue of
(\ref{limeq}) we multiply both parts of (\ref{limeq}) by $\eta_t$
in $H_0$. Integrating the resulting equality on $[0,\tau]$ and
observing that $\langle \left(\left(\bar{\mathcal W}_t+{\mathcal
H}_t\right)
    \cdot \nabla \right) \eta_t, \eta_t \rangle
    =0$ we conclude

 $$\frac{1}{2}\|\eta_\tau\|^2_0 +\nu \int_0^\tau \|\eta_t\|^2_1dt= \frac{1}{2}\|\eta_0\|^2_0+\int_0^\tau\langle \left({\mathcal H}_t
    \cdot \nabla \right) \bar{w}_t , \eta_t \rangle dt.$$

The summand $\langle \left({\mathcal H}_t
    \cdot \nabla \right) \bar{w}_t, \eta_t\rangle$ can be
    estimated as in the Section~\ref{cft}:
$$|\langle \left({\mathcal H}_t
    \cdot \nabla \right) \bar{w}_t, \eta_t\rangle| \leq
    C \|{\mathcal H}_t\|_1 \|\nabla \bar{w}_t\|_1 \|\eta_t\|_0 \leq C \|\nabla \bar{w}_t\|_1 \|\eta_t\|^2_0
    .$$

 One concludes  $\|\nabla \bar{w}_t\|_{1} \leq c'\|\bar{w}_t\|_{2}$,   while $\|\bar{w}_t\|_{2}$ are equibounded according to the
 Proposition~\ref{bound}.

Thus we get
$$ \frac{1}{2}\|\eta_\tau \|_0^2   \leq  \frac{1}{2}\|\eta_0
\|_0^2+c\int_0^\tau \|\eta_t\|^2_0
 dt,$$
and by application of Gronwall inequality we conclude
$$\frac{1}{2}\|\eta_\tau \|_0^2   \leq  \frac{1}{2}\|\eta_0
\|_0^2e^{cT}. \ \Box$$

\section{Proof of the Theorem~\ref{procon}}
\label{proc}

The proof of the result regarding $L_2$-approximate
controllability (Theorem~\ref{apcon}) for 2D NS system can be
found in \cite{AS43}; it holds also for 2D Euler system. Here we
provide a proof of the Theorem~\ref{procon}, which regards
controllability in finite-dimensional projection. Following the
steps of our proof the readers can  recover  the proof of
Theorem~\ref{apcon}.

Let $\mathcal L$ be a $\ell$-dimensional subspace of $H_2$ and
$\Pi^{\mathcal L}$ be $L_2$-orthogonal projection of $H_2$ onto
${\mathcal L}$. We start with constructing a finite-dimensional
coordinate subspace which is projected by $\Pi^{\mathcal L}$ onto
${\mathcal L}$.

To find one it suffices  to pick a $(\dim {\mathcal L}) \times
(\dim {\mathcal L})$-sub-matrix, from the $(\dim {\mathcal L})
\times \infty$ matrix which is a coordinate representation of
$\Pi^{\mathcal L}$. We look for more: for each $\eps >0$ we would
like to find a finite-dimensional coordinate subspace, which
contains an $\ell$-dimensional (non-ccordinate) subsubspace
${\mathcal L}_\eps$, which is $\eps$-close to ${\mathcal L}$. The
latter means that not only $\Pi^{\mathcal L}{\mathcal
L}_\eps={\mathcal L}$ but also $\Pi^{\mathcal L}|_{{\mathcal
L}_\eps}$ is $\eps$-close to the identity operator.

To achieve this we choose  an orthonormal basis $e_1, \ldots ,
e_\ell$ in ${\mathcal L}$ and takes for each $e_i$ its
finite-dimensional component (truncation) $\bar{e}_i$ which is
$\eps$-close to $e_i$. All $\bar{e}_i$ belong to some
finite-dimensional coordinate subspace $\mathcal S$ of $H_2$;
which reunites modes indexed by some symmetric set $S \subset
\mathbb{Z}^2$. Let $\Pi_S$ be $L_2$-orthogonal projection of $H_2$
onto $\mathcal S$. The subspace $\mathcal S$ together with the
subsubspace ${\mathcal L}_\eps$ spanned by $\bar{e}_1, \ldots ,
\bar{e}_\ell$ are the ones we looked for. Indeed
$$\|\Pi^{\mathcal L}\bar{e}_i-\bar{e}_i\|= \|\sum_{j=1}^\ell \langle \bar{e}_i, e_j\rangle
e_j-\bar{e}_i\|\leq \sum_{j=1}^\ell |\langle \bar{e}_i-e_i,
e_j\rangle| +\|e_i-\bar{e}_i\|\leq (\ell+1)\eps .$$

Without lack of generality we may assume that
$\|\Pi_S(\tilde{\varphi})-\tilde{\varphi}\|_0 \leq \eps$.

The set $\mathcal{K}^1$ of controlled modes is saturating, i.e.
for $\mathcal{K}^j$ defined by (\ref{setk2}), $\mathcal{K}^M
\supseteq S$ for some $M$. This means that the system is solidly
controllable in the observed component $q^S$.

In the proof of the Theo\-rem~\ref{fulcon} (Sec\-tion~\ref{cr}) we
started with a "full-di\-men\-sional" set of controlled modes
indexed by $\mathcal{K}^M$ and then constructed successively
controls which only enter the equations for the modes indexed by
${\mathcal K}^{M-1}, \ldots , {\mathcal K}^1$.

Assume that we are at the first induction step under the
conditions of the Lemma~\ref{m=1}, i.e. that all the coordinates
of the component $q^S$ are controlled. Following Lemma~\ref{m=1}
let us construct a family of controls which steers the
$q^S$-components of the corresponding trajectories $w_t$
 from $\Pi_S(\tilde{\varphi})$ to the points of the "ball"
$\mathcal{C}_R$ in $\mathcal S$. We can construct these controls
to actuate on an interval of arbitrarily small length $\tau >0$.
Denoting by $Q_.$ the component of $w_t$'s which is orthogonal to
$q^s$ we can conclude from (\ref{resc}) and (\ref{0resc}):
$$\|Q_t\|_0 \leq \|Q_0\|_0+C\tau , \ t \in [0,\tau], $$
for some constant $C>0$.

 Recall that $\|Q_0\|_0 =\|\Pi_S(\tilde{\varphi})-\tilde{\varphi}\|_0 \leq
\eps$.  Choosing $\tau \leq \eps/C$, we conclude $\|Q_t\|_0 \leq
2\eps, \ \forall t \in [0,\tau]$.

Let us check what happens with the component $Q_\cdot$ at generic
induction step of the proof of the Theorem~\ref{fam}. At the first
stage of each step (Subsection~\ref{stepgen}) we apply the
Approximation Lemma (Proposition~\ref{alem}). At this stage the
trajectories are approximated up to arbitrary small (uniformly for
$t \in [0,\tau]$) error $\delta >0$. We can choose $\delta \leq
\eps/(2M)$.

At the second stage of each induction step
(Subsection~\ref{stepfin}) the component $Q_\cdot$ (which belongs
to the image of the projection $\Pi_2$) suffers arbitrarily small
 alteration. We can make it (uniformly for  $t \in [0,\tau]$) smaller than
$\eps/(2M)$.

Therefore at each induction step the component $Q_\cdot$ suffers
alteration by value $\leq \eps/M$; total alteration is $\leq
\eps$. Hence after the induction procedure $\|Q_\tau\|_0 \leq
2\eps+\eps=3\eps$.

As a result we constructed a family of controls $x(\cdot;b)$ such
that the map $b \mapsto \Pi_S \circ \mathcal{F/T}_{T})(x(t;b))$
covers solidly the ball $\mathcal{C}_R$ in $\mathcal S$. Besides
$\|(\Pi^\perp_S \circ \mathcal{F/T}_{T})(x(t;b))\| \leq 3\eps$. If
$\eps >0$ is sufficiently small, then by construction of $S$ the
map $b \mapsto \Pi_{\mathcal{L}} \circ \Pi_S \circ
\mathcal{F/T}_{T})(x(t;b))$ covers the set $\mathcal{C}_{R/2} \cap
\mathcal{L}. \ \Box$

\section{Appendix}

\subsection{Non-classically forced
2D Euler equation: existence and uniqueness of so\-lu\-tions}
\label{exeul}
 We outline the proof which is a modification
 of the proof of the existence and uniqueness theorem for 2D Euler equation
 to be found in \cite{Ka}. Recall that the original proof of existence
 and uniqueness of classical solutions has been accomplished by
 W.Wolibner in \cite{Wo}.

Consider the nonclassically forced equation (\ref{ef}) with $\nu
=0$:
$$
\partial y_t/\partial t=\left({\mathcal Y}_t
    \cdot \nabla \right) y_t+\left({\mathcal V}^1_t
    \cdot \nabla \right) y_t+\left({\mathcal Y}_t
    \cdot \nabla \right) V^2_t+ V^0_t,$$
    where  $V^j_t \ (j=0,1,2)$ are trigonometric polynomials and
    ${\mathcal V}^1_t$ and ${\mathcal
Y}_t$ are divergence-free solutions
    of the equations
$$ \nabla^\perp \cdot {\mathcal V}^1_t=V^j_t,\ \nabla^\perp \cdot {\mathcal
Y}_t=y_t,$$  under  periodic boundary conditions.

Following the approach of \cite{Ka} let us introduce a map
$\xi_\cdot \mapsto \Phi(\xi)=\eta_\cdot $ which is defined by
means of the {\em linear} differential equations
\begin{eqnarray*}%%\label{map1}
\nabla^\perp \cdot \zeta_t=\xi_t, \\
 \partial \eta_t/\partial t=\left({\mathcal \zeta}_t
    \cdot \nabla \right) \eta_t+\left({\mathcal V}^1_t
    \cdot \nabla \right) \eta_t+\left({\mathcal \zeta}_t
    \cdot \nabla \right) V^2_t+ V^0_t.
%%\label{map2}
\end{eqnarray*}
It is easy to see that fixed points of the map $\Phi$ correspond
to classical solutions of the equation (\ref{ef}).

Choosing an appropriate set $\Omega$ of H\"olderian (of exponent
$\delta \in (0,1)$) with respect to time and space variables)
 functions with
$L^{(x)}_\infty$-norms bounded by a constant,  one is able to
establish, as in \cite{Ka}, that $\Phi$ maps $\Omega$ in itself.
Besides $S$ is compact convex subset of $C^0$ and existence of
fixed point is derived from Schauder theorem.

Analysis of the proof shows that the equiboundedness of the
$L^{(x)}_\infty$-norms of $V^j_t$ guarantee equiboundedness of the
$L^{(x)}_\infty$-norms of the corresponding solutions of
(\ref{ef}). This will prove the statement i) of the
Proposition~\ref{bound}.

\subsection{Forced
2D NS equation: existence, uniqueness and boundedness of
solutions}

The existence of  solutions from $L_\infty\left([0,t];H_2\right)$
for the nonclassically forced NS equation (\ref{ef}) can be
established in the same way as for classically forced NS equation,
for example by energy estimates for  Galerkin approximations.

In the same classical way we prove the boundedness of $\|y_t\|_2$,
and of $\int_0^T\|\frac{d}{dt}y_t\|_1^2dt$ i.e. the estimates
(\ref{eqbh1}) and (\ref{intdy1}). The boundedness of
$L^{(x)}_\infty$-norms (the estimate (\ref{eqblin})) follows then
from Sobolev inequality (see  \cite{Ada}).

\subsection{Proof of the Theo\-rem~\ref{bound}: equiboundedness of solutions for 2D Euler equation}

For the nonclassically forced 2D Euler equation the
$L_\infty^{(x)}$-equiboundedness of solutions
(Theorem~\ref{bound}; item i)) comes with the proof of existence
(see Subsection~\ref{exeul}).

To prove the statement ii) of the Theorem~\ref{bound}  we observe
first that uniform (in $t$) $L_\infty$-equi\-bo\-un\-ded\-ness of
$y_t$ implies their uniform (in $t$)
$H_0$-equi\-bo\-un\-ded\-ness. To arrive to the conclusion of the
assertion ii) let us differentiate both sides of the equation
(\ref{ef}), say, with respect to $x_i$. Abbreviating
$\partial/\partial x_i$ to $\partial_i$ we get:
\begin{eqnarray*}
\frac{\partial}{\partial t}(\partial_i y_t)=\left(\left({\mathcal
Y}_t
    +{\mathcal V}^1_t\right)
    \cdot \nabla \right) (\partial_i y_t)+\left(\left(\partial_i\left({\mathcal Y}_t
    +{\mathcal V}^1_t\right)\right)
    \cdot \nabla \right) y_t+  \\
    +\left(\left(\partial_i {\mathcal Y}_t\right)
    \cdot \nabla \right) V^2_t+ \left({\mathcal Y}_t
    \cdot \nabla \right) (\partial_i V^2_t) + \partial_i V^0_t.
\label{p1y}
\end{eqnarray*}
Multiplying both sides of the latter equality by $\partial_i y_t$
in $H_0$ we obtain
\begin{eqnarray}
\frac{1}{2}\frac{\partial}{\partial t} \|\partial_i y_t\|^2_0=
\nonumber \\ = \left\langle \left(\left({\mathcal Y}_t
    +{\mathcal V}^1_t\right)
    \cdot \nabla \right) \partial_i y_t, \partial_i y_t \rangle   +\langle \left(\left(\partial_i\left({\mathcal
    Y}_t+
    {\mathcal V}^1_t\right)\right)
    \cdot \nabla \right) y_t, \partial_i y_t \right\rangle
  \nonumber \\
    + \langle \left(\partial_i {\mathcal Y}_t
    \cdot \nabla \right) V^2_t, \partial_i y_t \rangle
    + \langle \left({\mathcal Y}_t
    \cdot \nabla \right) \partial_i V^2_t, \partial_i y_t \rangle +
    \langle \partial_i V^0_t, \partial_i y_t \rangle .
 \label{yp1y}
\end{eqnarray}

At the right-hand side of (\ref{yp1y}) the summand $\left\langle
\left(\left({\mathcal Y}_t
    +{\mathcal V}^1_t\right)
    \cdot \nabla \right) (\partial_i y_t), \partial_i y_t \right\rangle $
 is known to vanish, while the summand $\langle \left(\left(\partial_i\left({\mathcal Y}_t
    +{\mathcal V}^1_t\right)\right)
    \cdot \nabla \right) y_t, \partial_i y_t \rangle$ admits an upper estimate:

\begin{eqnarray}
\label{21}
    \left\langle \left(\left(\partial_i\left({\mathcal Y}_t
    +{\mathcal V}^1_t\right)\right)
    \cdot \nabla \right) y_t, \partial_i y_t \right\rangle \leq \nonumber \\ \leq C \left\|\partial_i\left({\mathcal Y}_t
    +{\mathcal V}^1_t\right)\right\|_{L_\infty}\|\nabla y_t\|_{L_2} \|\partial_i
    y_t\|_{L_2}\leq \\
    \leq C' \left\|\left(\partial_i\left({\mathcal Y}_t
    +{\mathcal V}^1_t\right)\right)\right\|_{L_\infty} \nonumber
    \|y_t\|^2_{H_1}.
\end{eqnarray}

Evidently $\|(\left(\partial_i{\mathcal Y}_t \right)\|_{L_\infty}
\leq c\|y_t\|_{L_\infty}$ and since, by virtue of i),
$\|y_t\|_{L_\infty}$ are
    bounded, then the upper estimate (\ref{21}) can be changed to
    $$\langle \left(\partial_i\left({\mathcal Y}_t
    +{\mathcal V}^1_t\right)
    \cdot \nabla \right) y_t, \partial_i y_t \rangle \leq C''\|y_t\|^2_{H_1}.  $$

The summand $\langle \left(\partial_i {\mathcal Y}_t
    \cdot \nabla \right) V^2_t, \partial_i y_t \rangle$ can be
    estimated from above by $$a\|\partial_i {\mathcal
    Y}_t\|_{L_2}\|\nabla  V^2_t\|_{L_\infty}\|\partial_i
    y_t\|_{L_2}.$$
     As long as $V^2_t$ is trigonometric polynomial in $x$ we can change the latter estimate to $
    a'\|y_t\|^2_{H_1}$. A similar upper estimate is valid for the
    summand $\langle \left({\mathcal Y}_t
    \cdot \nabla \right) (\partial_i V^2_t), \partial_i y_t
    \rangle$.

Finally $\langle \partial_i V^0_t, \partial_i y_t \rangle$ admits
an upper estimate
$$\alpha\left(\|\partial_i V^0_t\|_{L_2}^2+ \|\partial_i y_t\|_{L_2}^2 \right) \leq
\alpha'+\alpha''\|y_t\|^2_{H_1}.$$ Then we come to the
differential inequality  for $\|y_t\|^2_{H_1}$ denoted for brevity
by $\|y_t\|^2_1$:
$$\frac{\partial}{\partial t} \|
y_t\|^2_1 \leq  c'+c\| y_t\|^2_1,$$ wherefrom by the application
of the Gronwall inequality we conclude
$$\|y_t\|_1^2 \leq \|y_0\|_1^2e^{ct}+(c'/c)(e^{ct}-1),$$
and consequently $\sup_{t \in [0,T]}\|y_t\|_1 \leq b$ for some
$b>0$.

To arrive to the  estimate (\ref{eqbh1}) for $\|y_t\|_2$ (given
that the initial value $y_0$ belongs to $H_2$) we have to derivate
(\ref{p1y}) with respect to $x_j$, arriving to a differential
equation for $\partial_j\partial_i y$. Multiplying both parts of
this equation by $\partial_j\partial_i y$ we obtain
$$\frac{\partial}{\partial t}\|\partial_j \partial_i
y\|_0^2= \langle \left(\left({\mathcal Y}_t
    +{\mathcal V}^1_t\right)
    \cdot \nabla \right) \partial_j\partial_i y_t), \partial_j \partial_i y_t \rangle + \cdots .$$
The first term at the right-hand side vanishes and then the needed
estimate for $\|y_t\|_2^2$ is derived from the estimate for
$\|y_t\|_1^2$ by application of Young and Gronwall inequalities.

The integral estimate  iii) and even a stronger "pointwise"
estimate for $\left\|\frac{\partial}{\partial t}(\partial_i
y_t)\right\|$ can be concluded   from (\ref{p1y}). Indeed
\begin{eqnarray*}
\left\|\frac{\partial}{\partial t}(\partial_i y_t)\right\| \leq
C\left(\|{\mathcal Y}_t
    +{\mathcal V}^1_t\|_{L_\infty}\|\nabla \partial_i y_t  \|+ \|\partial_i\left({\mathcal Y}_t
    +{\mathcal V}^1_t\right)\|_{L_\infty}\|\nabla  y_t  \| \right)+ \\
    +C_1\|y_t\|_1+c_2 \leq C'(1+\|y_t\|_0)\|y_t\|_2+c'_2.
\end{eqnarray*}

%%%%%%%%%%%%%%%%%%%%%%%%%%%%%%%%%%%%%%%%%%%%%%%%%%%%%%%%%%%%%%%%%%%%

\subsection{Continuity with respect to relaxation metric: proof of the Theo\-rem~\ref{corr}}

Pick an element $(\bar{V}^0_t,\bar{V}^1_t,\bar{V}^2_t)$ from
$\mathbf{F}_B$ and denote by $\by_t$ the solution of the equation
\begin{equation}\label{2kb}
  \partial_t \by_t =\left(\bar{\mathcal Y}_t \cdot \nabla \right) \by_t
   + \left(\bar{{\mathcal V}}^1_t
    \cdot \nabla \right) \by_t+\left(\bar{{\mathcal Y}}_t
    \cdot \nabla \right) \bar{V}^2_t+ \nu \Delta \by_t+\bar{V}^0_t.
\end{equation}

  Let $y_t$ be a solution of the "perturbed" equation
\begin{eqnarray}\label{delv}
\partial_t y_t =\left({\mathcal Y}_t \cdot \nabla \right) y_t
   + \left(\left(\bar{{\mathcal V}}^1_t+{\mathcal V}^1_t\right)
    \cdot \nabla \right) y_t+\left({\mathcal Y}_t
    \cdot \nabla \right) \left(\bar{V}^2_t+V^2_t\right)+ \\ +\nu
    \Delta y_t+
    \bar{V}^0_t+V^0_t. \nonumber
\end{eqnarray}
Recall that $\bar{{\mathcal V}}^1_t$ is the divergence-free
solution of the equation: $\nabla^\perp  \bar{{\mathcal
V}}^1_t=\bar{\mathcal V}^1_t$.

Subtracting (\ref{2kb}) from (\ref{delv}) and introducing the
notation
$$\eta_t=y_t-\bar{y}_t, \ {\mathcal H}_t={\mathcal Y}_t
-\bar{\mathcal Y}_t,$$  we obtain the equation for $\eta_t$:
\begin{eqnarray}\label{eqeta}
\partial_t \eta_t
    =\left({\mathcal Y}_t
    \cdot \nabla \right) \eta_t+\left({\mathcal H}_t \cdot \nabla \right) \by_t+
\left(\left(\bar{{\mathcal V}}^1_t+{\mathcal V}^1_t\right)
    \cdot \nabla \right) \eta_t+\left({\mathcal V}^1_t
    \cdot \nabla \right)\by_t+ \nonumber \\
    +\left({\mathcal Y}_t \cdot \nabla
    \right)
    V^2_t+   \left({\mathcal H}_t
    \cdot \nabla \right)\bar{V}^2_t +\nu \Delta \eta_t+V^0_t.
\end{eqnarray}

We would like to evaluate $\|\eta\|_0$; to this end we multiply in
$H_0$ both sides of (\ref{eqeta}) by $\eta_t$. At the right-hand
side we obtain $\frac{1}{2} \partial_t\|\eta_t\|_0^2$, while at
the right-hand side  the terms $\langle \left({\mathcal Y}_t
    \cdot \nabla \right) \eta_t, \eta_t \rangle$ and $\langle \left(\left(\bar{{\mathcal V}}^1_t+{\mathcal V}^1_t\right)
    \cdot \nabla \right) \eta_t, \eta_t \rangle$ both vanish. Taking into account
    that $\langle \Delta \eta_t, \eta_t \rangle \leq 0$ at the right-hand side,  we arrive
    to the inequality:
\begin{eqnarray*}
\frac{1}{2}\partial_t \|\eta_t\|_0^2 \leq   \langle
\left({\mathcal H}_t \cdot \nabla \right)
\left(\by_t+\bar{V}^2_t\right), \eta_t \rangle+ \\  + \langle
\left({\mathcal V}^1_t    \cdot \nabla \right)\by_t, \eta_t
\rangle +\langle \left({\mathcal Y}_t \cdot \nabla    \right)
    V^2_t, \eta_t \rangle  + \langle
    V^0_t, \eta_t \rangle .
\end{eqnarray*}
Hence
\begin{eqnarray}\label{eta0}
\frac{1}{2}\|\eta_\tau\|_0^2 \leq  \frac{1}{2}\|\eta_0\|_0^2+
\int_0^\tau \langle \left({\mathcal H}_t \cdot
\nabla \right) \left(\by_t+\bar{V}^2_t\right), \eta_t \rangle dt+ \\
+ \int_0^\tau \left(\langle \left({\mathcal V}^1_t    \cdot \nabla
\right)\by_t, \eta_t \rangle +\langle \left({\mathcal Y}_t \cdot
\nabla \right)
    V^2_t, \eta_t \rangle  + \langle
    V^0_t, \eta_t \rangle\right) dt . \nonumber
\end{eqnarray}

  What for the first integrand  in the right-hand side, then
\begin{eqnarray*}
\left|\langle \left({\mathcal H}_t \cdot \nabla \right)
\left(\by_t+\bar{V}^2_t\right), \eta_t \rangle \right| \leq \\
\leq  c\|{\mathcal H}_t\|_1 \|\nabla
\left(\by_t+\bar{V}^2_t\right)\|_1\|\eta_t\|_0 \leq
c'\|\bar{y}_t+\bar{V}^2_t\|_2\|\eta_t\|^2_0.
\end{eqnarray*}

As long as $\bar{V}^2_t$ are trigonometric polynomials with
uniformly bounded coefficients and according to
Theorem~\ref{bound} $\|\bar{y}_t\|_2$ are equibounded, then the
latter estimate can be changed to $c''\|\eta_t\|^2_0$.

All terms of the second integrand at the right-hand side of
(\ref{eta0}) contain "factors" $V^0_t,V^1_t,V^2_t$ which are small
in relaxation metric. To estimate this integral one can use
Lemma~\ref{ss:105}. Its assumptions are verified as far as the
values
$$\int_0^T\left\|\frac{d}{dt}\partial_i (\by_t)_j\right\|_0^2dt, \int_0^T\left\|\dot{\eta}_t\right\|_0^2dt,
\int_0^T\left\|\dot{y}_t\right\|_0^2dt$$
 are equibounded. Say, the value of the integral $\int_0^\tau \left(\langle \left({\mathcal V}^1_t    \cdot \nabla
\right)\by_t, \eta_t \right\rangle dt$ is small, because
${\mathcal V}^1_t$ is small in relaxation norm, and
$\int_0^T\left(\frac{d}{dt}\left(\partial_i (\by_t)_j
(\eta_t)_j\right)\right)^2dt$ are bounded.

For any $\delta >0$ we can take $V^0_t,V^1_t,V^2_t$ sufficiently
small in relaxation metric, in such a way that  (\ref{eta0})
implies
$$\frac{1}{2}\|\eta_\tau\|_0^2 \leq  \frac{1}{2}\|\eta_0\|_0^2+
c''\int_0^\tau \|\eta_t\|_0^2 dt+\delta .$$ Then the smallness of
$\|\eta_t\|_0$ is concluded by application of the Gronwall
inequality. $\Box$

% ----------------------------------------------------------------

\end{document}